\documentclass[11pt]{article}
\usepackage{amssymb, latexsym, amsmath}
\usepackage{amsfonts, graphics}
\def\R{\mathbb R}

\begin{document}
\sloppy
\newtheorem{theorem}{Theorem}[section]
\newtheorem{definition}[theorem]{Definition}
\newtheorem{proposition}[theorem]{Proposition}
\newtheorem{example}[theorem]{Example}
\newtheorem{remark}[theorem]{Remark}
\newtheorem{cor}[theorem]{Corollary}
\newtheorem{lemma}[theorem]{Lemma}

\renewcommand{\theequation}{\arabic{section}.\arabic{equation}}

\title{Steady states of the spherically symmetric Vlasov-Poisson system
as fixed points of a mass-preserving algorithm}

%\title{Steady states of the spherically symmetric Vlasov-Poisson system
%obtained as fixed points of a mass-preserving algorithm}

%\title{A new method to prove the existence 
%\\ of static solutions to the spherically symmetric Vlasov-Poisson system} 

\author{H{\aa}kan Andr\'{e}asson\\
        Mathematical Sciences\\
        Chalmers University of Technology\\
        University of Gothenburg\\
        SE-41296 Gothenburg, Sweden\\
        email: hand@chalmers.se\\
        \ \\
        Markus Kunze\\
        Mathematisches Institut\\
        Universit\"at K\"oln\\
        Weyertal 86-90\\
        D-50931 K\"oln, Germany\\
        email: mkunze1@uni-koeln.de\\        
        \ \\
        Gerhard Rein\\
        Fakult\"at f\"ur Mathematik, Physik und Informatik\\
        Universit\"at Bayreuth\\
        D-95440 Bayreuth, Germany\\
        email: gerhard.rein@uni-bayreuth.de}

\maketitle

\begin{abstract} 
We give a new proof for the existence of spherically symmetric steady states 
to the Vlasov-Poisson system, following a strategy that has been used
successfully 
to approximate axially symmetric solutions numerically,
both to the Vlasov-Poisson system 
and to the Einstein-Vlasov system. There are several reasons why a mathematical analysis 
of this numerical scheme is important. A generalization of the present result to the case 
of flat axially symmetric solutions would prove that the steady states obtained numerically 
in \cite{AR3} do exist. Moreover, in the relativistic case the question
whether a steady state can be obtained by this scheme seems to be related
to its dynamical stability. 
This motivates the desire for a deeper understanding of this strategy. 
\end{abstract}

\section{Introduction}

\setcounter{equation}{0}

Spherically symmetric static solutions and axially symmetric stationary solutions 
of the Vlasov-Poisson system and of the Einstein-Vlasov system have been extensively studied, 
both analytically and numerically. A key difference exists between the approaches used 
to obtain solutions in the spherically symmetric case versus the axially symmetric case. 
In the spherically symmetric setup, the equations can be treated as an initial value problem. 
Starting with data at $r=0$, a nonlinear system of integro-differential equations 
can be solved directly, at least numerically, allowing for the construction 
of a wide range of solutions; cf.~\cite{An0,AFT,AR2,AR3,RaRe,R1}. 
In contrast, the axially symmetric case requires a different approach. Numerical methods employed in previous works rely on iterative procedures that must converge to produce a solution. 

A crucial ingredient for numerical convergence is that the mass of the
iterates  
is kept constant. We call this a mass-preserving algorithm. This algorithm has successfully 
been used to construct stationary solutions numerically in the axially symmetric case 
in a number of works, such as \cite{AAL1,AAL2,AR1,ST1,ST2,ST3}.
An important feature in these works 
is that the solutions that are obtained are far from being
spherically symmetric. This is very different for the axially symmetric
solutions whose existence has been proven, since these results rely
on the implicit function theorem, perturbing off a spherically symmetric solution \cite{AKR1,AKR2,R2,S}, and hence the resulting
solutions are close to being spherically symmetric. 
In the case of the Vlasov-Poisson system, also flat axially symmetric steady states have been studied. 
The mass-preserving algorithm was used in \cite{AR3} to numerically construct steady states 
with the property that the corresponding rotation curves are flat, 
in agreement with observations for disk galaxies. Existence results 
for flat axially symmetric steady states have been obtained in \cite{FR,R2}, 
but these results do not cover the type of solutions constructed numerically
in \cite{AR3}. 

The purpose of this work is to analyze the mass-preserving algorithm 
that has shown to be so powerful for constructing stationary solutions numerically. 
This algorithm can be formulated as a fixed point problem,
and the main achievement 
of the present paper is the proof that
under suitable assumptions a fixed point to this problem exists 
in the case of the spherically symmetric Vlasov-Poisson system.
We have therefore found a new method to analytically construct steady states.

In view of the discussion above there are more important reasons
for the present investigation. 
We believe that our result can be extended to cover 
the flat axially symmetric case, which would prove that the solutions 
obtained numerically in \cite{AR3} do actually exist. 
The estimates in the present work are tailored to the spherically
symmetric case, 
and the modifications required for the flat case
are not straightforward and pose an interesting open problem. 
Another desirable generalization of our method concerns the spherically symmetric Einstein-Vlasov 
system.
Although the existence of solutions can be obtained by other means,
there is numerical evidence of an intriguing relation between the fact
that a steady state can be obtained by the mass-preserving algorithm 
and its dynamical stability.
Indeed, both stable and unstable steady states are studied in \cite{AR1} 
and one can try to numerically construct these steady states by applying
the mass-preserving algorithm. 
There is evidence that only in the case of stable steady
states the algorithm converges; 
see also the discussion in \cite{AA}.
Analytic progress on the mass-preserving algorithm 
in the relativistic case could shed some
light on this observation, and it could even be valuable 
for making progress on the nonlinear stability problem. 

The outline of the paper is as follows.
In the next section the fixed point problem 
for the Vlasov-Poisson system is formulated in quite general terms.   
In Section \ref{sect_T} the assumption of spherical symmetry is added. 
In Section \ref{sect_ubcc} a crucial bound on the potential is derived. 
Finally, in Section \ref{sect_compT} 
the existence of a fixed point is proven and the main result,
Theorem \ref{main_thm}, is established. 

\section{The setup}

Our aim is to prove the existence of spherically symmetric or axially
symmetric, possibly flat steady states of the Vlasov-Poisson
system by the following general strategy.
We prescribe an ansatz for the particle distribution
\begin{equation} \label{ansatz}
   f = \Phi(E,L)
\end{equation}
in terms of the particle energy
\[ E = \frac{1}{2} |v|^2 + U(x) \]
and
\[ L = |x\times v|^2 \ \mbox{or}\ L = x_1 v_2 - x_2 v_1, \]
the modulus squared of angular momentum or its third component, respectively.
For a given spatial density $\rho = \rho(x)$, which is spherically 
or axially symmetric and has finite mass $M>0$ and compact support, 
we define its induced gravitational potential as
\[ U_\rho(x) = -\int \frac{\rho(y)}{|x-y|} dy. \]
We substitute $U_\rho$ into the ansatz \eqref{ansatz} to
obtain a new spatial density
\[
\tilde \rho (x) = \int \Phi\left(\frac{1}{2} |v|^2 + U_\rho (x), L\right)\,dv.
\]
In addition, we define an amplitude
\[
A(\tilde \rho) = M \left(\int \tilde \rho(x)\, dx \right)^{-1},
\]
so that the new spatial density $A(\tilde \rho)\, \tilde \rho$ again
has mass $M$.
Now assume that we have a fixed point $\rho^\ast$ of the following map,
defined on some suitable domain $D$:
\[
T \colon D \ni \rho \mapsto U_\rho \mapsto
   \tilde \rho = \int \Phi\left(\frac{1}{2} |v|^2 + U_\rho, L\right)\, dv
   \mapsto A(\tilde \rho)\, \tilde \rho \in D.
\]
Then clearly
\[ \rho^\ast= A(\rho^\ast)^{-1}
   \int \Phi\left(\frac{1}{2} |v|^2 + U_{\rho^\ast}, L\right)\, dv \]
so that $(f^\ast,\rho^\ast,U_{\rho^\ast})$ is a steady state of
the Vlasov-Poisson system, where $f^\ast$ is given by the new ansatz
\[ f^\ast = A(\rho^\ast)^{-1} \Phi(E,L). \]
The above observation is the basis for several numerical schemes that are used 
to construct axially symmetric solutions, as described in the introduction. 
We will show that it can also employed to prove the
existence of steady states, at least in the spherically symmetric case.

\section{The operator $T$ in the spherically symmetric case}
\label{sect_T}

\setcounter{equation}{0}

As a proof of concept for the approach described above, 
we consider spherically symmetric, isotropic steady states, 
where we restrict ourselves to the polytropic ansatz
\begin{equation} \label{polyansatz}
   f(x,v) = (E_0 - E)_+^k, 
\end{equation}
for $E_0<0$ denoting some fixed cut-off energy and $k>-1$, 
which will be further restricted in what follows. 
We wish to obtain steady states of fixed, prescribed
mass $M>0$ supported in a fixed ball of radius $R>0$. As explained above,
the resulting steady state will actually have a microscopic equation
of state where the right hand side in \eqref{polyansatz}
is multiplied by some positive constant. To simplify notation
we will restrict ourselves to the choices
\[ M=1,\ R=1,\ E_0=-1. \]
Letting $B=B_1(0)$ denote the open ball in $\R^3$
of radius $1$ centered at the origin, we introduce the following set: 
\begin{eqnarray}\label{space}
  D := \{ \rho\in L^p(B)
  &\mid&
  \rho\geq 0\ \mbox{is spherically symmetric
    and decreasing with} \nonumber\\
  && \mbox{$\int_B \rho\,dx=1$}\}.
\end{eqnarray}
Here $p>3/2$, and spherical symmetry means that $\rho(x)=\rho(r)$
where $r=|x|$; in what follows we will continue to use the
identification encoded in this abuse of notation. Moreover,
we will think of functions in $D$ to be extended by $0$ to the whole
space $\R^3$.
The set $D$ is a closed, convex subset of the Banach space $L^p(B)$.
For $\rho\in D$,
\[ U_\rho(x) = -\int \frac{\rho(y)}{|x-y|} dy
   = -\frac{4\pi}{r} \int_0^r \rho(s)\, s^2 ds
   - 4\pi \int_r^\infty \rho(s)\, s \, ds, \]
and as a function of $r$, $U_\rho$ is increasing. 
Moreover, $U_\rho (1) = -1$.
Hence $U_\rho$ is the unique solution to the boundary value problem
\[ \Delta U = 4 \pi \rho\ \mbox{on}\ B,\ U=-1\ \mbox{on}\ \partial B. \]
Since $\rho\in L^p(B)$ it follows that $U_\rho\in W^{2,p}(B)$, 
and due to $p>3/2$ the latter space is compactly embedded into the space
$C(\overline{B})$. Hence the mapping
$D \ni \rho \mapsto U_\rho \in C(\overline{B})$ is a compact operator.

The function
\[ \tilde \rho(x) = \int \left(-1-\frac{1}{2}|v|^2 -U_\rho(x)\right)_+^k dv
   = c_k \left(-1 - U_\rho(x)\right)_+^{k+3/2}, \]
where
\[ c_k := 4\pi \sqrt{2} \int_0^1 \eta^k \sqrt{1-\eta}\, d\eta, \]
is continuous, non-negative, spherically symmetric, and decreasing,
and it vanishes outside $\overline{B}$ since $U_\rho(x)>-1$ there.
Moreover, the mapping
$C_b (\overline B)\ni U_\rho \mapsto \tilde \rho \in L^p(B)$
is continuous.
If we multiply $\tilde \rho$ by
\[ A(\tilde \rho) = \left(\int \tilde\rho(x)\, dx\right)^{-1} =
   \left(c_k \int_B\left(-1 - U_\rho(x)\right)_+^{k+3/2} dx\right)^{-1} \]
we obtain a density that is again contained in the set $D$,
and we need to control the behavior of this amplitude.
This will be done in the next section. 

\section{A uniform bound close to the center}
\label{sect_ubcc} 

\setcounter{equation}{0} 

\begin{lemma}\label{boundatcenter}
There exists $r_0\in ]0, 1]$ and $\delta>0$ such that 
for every $\rho \in D$ the induced potential has the property 
that $U_\rho (r_0)\le -1 -\delta$. In particular, 
$U_\rho (r)\le -1 -\delta$ for all $r\in [0, r_0]$. 
\end{lemma}

\noindent
{\bf Proof\,:}
For fixed $r_0\in ]0, \frac{1}{10}]$ we assert that 
\[ \delta=\min\Big\{\frac{1}{10}, \frac{1}{2}\Big(\frac{1}{r_1}-1\Big)\Big\}>0 \] 
is a suitable choice, where $r_1\in ]r_0, 1[$ satisfies
$r_1^3\ge\frac{4}{5}$. 
If $4\pi\int_0^{r_0} s^2\rho(s)\,ds\ge\frac{11}{10}\,r_0$, then 
\begin{eqnarray*} 
  U_\rho(r_0)
  & = &
  -\frac{4\pi}{r_0}\int_0^{r_0} s^2\rho(s)\,ds-4\pi\int_{r_0}^1 s\rho(s)\,ds\\
  & \le &
  -\frac{4\pi}{r_0}\int_0^{r_0} s^2\rho(s)\,ds\\
  & \le &
  -\frac{11}{10}=-1 -\frac{1}{10}\le -1 -\delta
\end{eqnarray*} 
as desired. Thus in what follows we can assume that 
\begin{equation}\label{klein}
   4\pi\int_0^{r_0} s^2\rho(s)\,ds < \frac{11}{10}\,r_0. 
\end{equation} 
Then $1-r_1^3\le\frac{1}{5}$ and hence
\begin{equation}\label{novap} 
  \frac{3}{1-r_1^3}\,\Big(\frac{1}{2}-\frac{11}{10}\,r_0\Big)\,
  \frac{r_1^3-r_0^3}{3}
  \ge 5\,\Big(\frac{1}{2}-\frac{11}{100}\Big)\,
  \Big(\frac{4}{5}-\frac{1}{1000}\Big)
  \ge\frac{3}{2}.
\end{equation} 
This implies that
\begin{equation}\label{lgtea} 
   4\pi\int_{r_0}^{r_1} s^2\rho(s)\,ds \geq \frac{1}{2}.
\end{equation} 
Indeed, suppose that 
$4\pi\int_{r_0}^{r_1} s^2\rho(s)\,ds < \frac{1}{2}$. 
Then $4\pi\int_0^1 r^2\rho(r)\,dr=1$ together with (\ref{klein}) yields
\begin{eqnarray*} 
  \frac{1}{2}
  & < &
  4\pi\int_0^{r_0} s^2\rho(s)\,ds + 4\pi\int_{r_1}^1 s^2\rho(s)\,ds \\
  & < &
  \frac{11}{10}\,r_0+\frac{4\pi}{3}\,\rho(r_1)\,(1-r_1^3).
\end{eqnarray*}
In the last estimate we used the fact that $\rho \in D$ by definition of $D$
implies that $\rho$ is decreasing, more precisely, there
exists a pointwise defined representative of $\rho$, which we continue to
call $\rho$, and which has the property that $\rho$ is decreasing on the
set $[0,1]\setminus N$ where $N$ is a set of measure zero. By making $r_1$
slightly bigger if necessary we can make sure that $r_1 \notin N$.
Hence
\[ \rho(r_1) > \frac{3}{4\pi (1-r_1^3)}\Big(\frac{1}{2}-\frac{11}{10}\,r_0\Big), \] 
from which we deduce that 
\[ \rho(r)\ge\frac{3}{4\pi (1-r_1^3)}\Big(\frac{1}{2}-\frac{11}{10}\,r_0\Big),
   \quad \mbox{for almost all}\ r\in [0, r_1]. \] 
As a consequence, by (\ref{novap}),  
\[
   \frac{1}{2}  >  4\pi\int_{r_0}^{r_1} s^2\rho(s)\,ds
    \ge  \frac{3}{1-r_1^3}\Big(\frac{1}{2}-\frac{11}{10}\,r_0\Big)
   \,\frac{r_1^3-r_0^3}{3}\ge\frac{3}{2},
\]
and this contradiction completes the proof of (\ref{lgtea}). 

Let $\alpha=4\pi\int_{r_0}^{r_1} s^2\rho(s)\,ds$. 
Then $\alpha\ge\frac{1}{2}$ by (\ref{lgtea}), and
\begin{eqnarray*}
  4\pi\int_{r_0}^1 s\rho(s)\,ds
  & = &
  4\pi\int_{r_0}^{r_1}\frac{s^2}{s}\,\rho(s)\,ds
  +4\pi\int_{r_1}^1 s\rho(s)\,ds\\
  & \ge &
  \frac{4\pi}{r_1}\int_{r_0}^{r_1} s^2\rho(s)\,ds
  +4\pi\int_{r_1}^1 s^2\rho(s)\,ds\\
  & = &
  \frac{\alpha}{r_1}+4\pi\int_{r_0}^1 s^2\rho(s)\,ds
  -4\pi\int_{r_0}^{r_1} s^2\rho(s)\,ds\\
  & = &
  \Big(\frac{1}{r_1}-1\Big)\alpha+4\pi\int_{r_0}^1 s^2\rho(s)\,ds\\
  & \ge &
  \frac{1}{2}\Big(\frac{1}{r_1}-1\Big)+4\pi\int_{r_0}^1 s^2\rho(s)\,ds\\
  & \ge &
  \bigg[\frac{1}{2}\Big(\frac{1}{r_1}-1\Big)+1\bigg]
  \,4\pi\int_{r_0}^1 s^2\rho(s)\,ds \\
  & \ge &
  (1+\delta)\,4\pi\int_{r_0}^1 s^2\rho(s)\,ds. 
\end{eqnarray*}
Hence 
\begin{eqnarray*} 
  U_\rho(r_0)
  & = &
  -\frac{4\pi}{r_0}\int_0^{r_0} s^2\rho(s)\,ds-4\pi\int_{r_0}^1 s\rho(s)\,ds \\
  & \le &
  -\frac{4\pi}{r_0}\int_0^{r_0} s^2\rho(s)\,ds
  -(1+\delta)\,4\pi\int_{r_0}^1 s^2\rho(s)\,ds \\
  & = &
  -\frac{4\pi}{r_0}\int_0^{r_0} s^2\rho(s)\,ds
  -(1+\delta)\,\bigg(1-4\pi\int_0^{r_0} s^2\rho(s)\,ds\bigg) \\
  & = &
  -(1+\delta)+4\pi\,\Big(1+\delta-\frac{1}{r_0}\Big)
  \int_0^{r_0} s^2\rho(s)\,ds.
\end{eqnarray*}
Since $r_0(1+\delta)\le\frac{1}{10}\,\frac{11}{10}=\frac{11}{100}\le 1$, 
we have $1+\delta-\frac{1}{r_0}\le 0$,
so that $U_\rho(r_0)\le -1-\delta$, 
which completes the proof. 
{\hfill$\Box$}\bigskip

\section{Compactness of $T$ and existence of a fixed point}
\label{sect_compT} 

\setcounter{equation}{0}

We continue the analysis of the operator $T$.
Using Lemma \ref{boundatcenter} it follows that
\[ \int_B\left(-1 - U_\rho(x)\right)_+^{k+3/2} dx
   \geq \frac{4\pi}{3} r_0^3 \delta^{k+3/2}. \]
On the other hand,
\[ U_\rho' (r) = \frac{4\pi}{r^2}\int_0^r\rho(s)\, s^2 ds \leq \frac{1}{r^2}, \]
and since $U_\rho(1)=-1$ we deduce that
\[ U_\rho(r) \geq - \frac{1}{r}. \]
This implies that
\[ \int_B\left(-1 - U_\rho(x)\right)_+^{k+3/2} dx
   \leq 4\pi\int_0^1 s^{2-k-3/2} ds < \infty, \]
provided that $k<3/2$. Thus there exist constants $C_1, C_2 >0$ such that for
every $\rho\in D$ the estimate
\begin{equation}\label{ampest}
   C_1 \leq A(\tilde \rho) \leq C_2
\end{equation}
holds. The mapping $D\ni \rho \mapsto A(\tilde \rho)$ 
is continuous, and by the previous estimate it is compact. 
Together with the results in Section \ref{sect_T} this implies that the operator
$T\colon D\to D$ is compact; it is continuous, and for any bounded sequence
$(\rho_n)\subset D$ the sequence $(T(\rho_n))$ has a convergent subsequence.
\smallskip

A further useful observation is as follows. 

\begin{lemma}
  Let $-1<k<1/2$. Then the range $T(D)\subset L^p(B)$ is bounded. 
\end{lemma} 
\noindent
{\bf Proof\,:} Recalling \eqref{ampest}, we only need to estimate ${\|\tilde \rho\|}_p$
uniformly on $D$. For $\rho \in D$ we know that $U_\rho(r)\geq -1/r$, and
this implies that
\[ {\|\tilde \rho\|}_p^p \leq C \int_0^1 r^{2-p(k+3/2)} dr < \infty, \]
provided
\[ p<\frac{3}{k+3/2}. \]
Since we required that $p>3/2$, this restricts the possible exponents $k$
of our ansatz to $-1<k<1/2$. {\hfill$\Box$}\bigskip 

Combining the previous results, we are able 
to state the main result of this paper. 

\begin{theorem}\label{main_thm} 
Let $-1<k<1/2$. Then the operator $T: D\to D$ has a fixed point. 
In particular, the Vlasov-Poisson system has a steady state 
for the ansatz function (\ref{polyansatz}). 
\end{theorem} 
\noindent
{\bf Proof\,:} We have seen above that $D\subset L^p(B)$ is closed and convex, $T: D\to D$ 
is compact and $T(D)\subset L^p(B)$ is bounded. Let $r=\sup\,\{{\|T(\rho)\|}_p: \rho\in D\}$ 
and $D_r=D\cap\{\rho\in L^p(B): {\|\rho\|}_p\le r\}$. Then $D_r\subset L^p(B)$ 
is closed, bounded and convex, and $T: D_r\to D_r$ is compact. 
Hence it has a fixed point by Darbo's Theorem; see \cite[Theorem 9.1]{deim}. 
{\hfill$\Box$}

\end{document}